\documentclass[11pt]{article}
\usepackage[utf8]{inputenc}
\usepackage[fleqn]{amsmath}
\usepackage{amsthm}
\usepackage{amssymb}
\usepackage[T1]{fontenc}
\usepackage{lmodern}
\usepackage[english]{babel}
\usepackage{extsizes}
\usepackage{xcolor}
\usepackage{fancybox}
\usepackage{graphicx}
\usepackage{sectsty}
\usepackage{enumerate}
\usepackage[french]{varioref}
\usepackage{mathtools}
\usepackage{mathrsfs}
\usepackage{blkarray}
\usepackage{framed}
\usepackage[numbers]{natbib}
\usepackage{mathbbol}
\usepackage{etoolbox}
\usepackage{wasysym}
\usepackage{caption}

\usepackage[ocgcolorlinks,colorlinks=true,urlcolor=blue]{hyperref}
\makeatletter
\pretocmd{\NAT@citexnum}{\@ifnum{\NAT@ctype>\z@}{\let\NAT@hyper@\relax}{}}{}{}
\makeatother

\newtheorem{Th}{Theorem}
\newtheorem{Prop}[Th]{Proposition}
\newtheorem{Le}[Th]{Lemma}
\newtheorem{Cor}[Th]{Corollary}
\theoremstyle{definition}

\newtheorem{Rk}{Remark}
\newtheorem{Ex}{Example}
\mathchardef\ordinarycolon\mathcode`\:
\mathcode`\:=\string"8000
\begingroup \catcode`\:=\active
  \gdef:{\mathrel{\mathop\ordinarycolon}}
\endgroup

\newcommand{\ind}[1]{\mathbb{1}_{#1}\,}

\newcommand{\RR}{\mathbf{R}}
\newcommand{\ZZ}{\mathbf{Z}}
\newcommand{\NN}{\mathbf{N}}

\newcommand{\conv}[2][n]{\underset{#1\rightarrow #2}{\longrightarrow}}

\newcommand{\EEE}[1]{\operatorname{\mathbb{E}}\left[\,#1\,\right]}

\newcommand{\PPP}[1]{\operatorname{\mathbb{P}}\left(#1\right)}
\newcommand{\PP}{\operatorname{\mathbb{P}}}

\newenvironment{prooft}[1]{\vskip 2mm\noindent {\bf Proof of #1.}} {\hfill
                    $\square$ \vskip 2mm \noindent}
\usepackage[top=0.85in, bottom=0.85in, left=1.20in, right=1.20in]{geometry}

\begin{document}
\pagenumbering{arabic}
\pagestyle{plain}
\date{}

\title{On the discrepancy of powers of random variables}
\author{Nicolas Chenavier\footnote{Universit\'e Littoral C\^ote d'Opale, EA 2797, LMPA, 50 rue Ferdinand Buisson, F-62228 Calais, France. E-mail: nicolas.chenavier@univ-littoral.fr, \textit{corresponding author}}, Dominique Schneider \footnote{Universit\'e Littoral C\^ote d'Opale, EA 2797, LMPA, 50 rue Ferdinand Buisson, F-62228 Calais, France. E-mail: dominique.schneider@univ-littoral.fr}}
\maketitle

\begin{abstract}
Let $(d_n)$ be a sequence of positive numbers and let $(X_n)$ be a sequence of positive independent random variables. We provide an upper bound for the deviation between the distribution of the mantissaes of $(X_n^{d_n})$  and the Benford's law. If $d_n$ goes to infinity at a rate at most polynomial, this deviation converges a.s. to 0 as $N$ goes to infinity. 
\end{abstract}





\section{Introduction}
A sequence of positive numbers $(x_n)$ is said to satisfy the first digit phenomenon if 
\[\lim_{N\rightarrow\infty}\frac{1}{N}\sum_{n=1}^N\ind{F(x_n)=k} = \log_{10}\left(1+\frac{1}{k} \right), k=1,\ldots, 9,\] where $F(x_n)$ is the first digit of  $x_n$, and where $\ind{A}$ denotes the indicator function of any subset $A$. Such a phenomenon was observed by Benford and Newcomb on real life numbers \cite{benford, newcomb}. It is extensively used in various domains, such as fraud detection \cite{nigrini}, computer design \cite{hamming} and image processing \cite{XWLD}. As an extension of the first digit phenomenon, the notion of Benford sequence is introduced as follows. Let $\mu_{10}$ be the measure on the interval $[1,10)$ defined by \[
\mu_{10}([1,a))= \log_{10} a, (1\leq a <10),
\]
where $\log_{10}a$ denotes the logarithm in base $10$ of $a$. Let ${\mathcal{M}}_{10}(x)$ be the mantissa in base $10$ of a positive number $x$, i.e. ${\mathcal{M}}_{10}(x)$  is the unique number  in $[1,10)\,$ such that there exists an integer $k$ satisfying $x={\mathcal{M}}_{10}(x)10^k$. A set of numbers $(x_n)$ is referred to as a Benford sequence if for any $1\leq a < 10$, we have \[\lim_{N\rightarrow\infty}\frac{1}{N}\sum_{n=1}^N\ind{\mathcal{M}_{10}(x_n)\in [1,a)} = \mu_{10}([1,a)).\] In particular, each Benford sequence satisfies the first digit phenomenon since $F(x)=k$ if and only if $\mathcal{M}_{10}(x)\in [k, k+1)$, with $x>0$, $k=1,\ldots, 9$. For instance, the sequences $(2^n)$, $(n!)$ and $(n^n)$ are Benford. For various examples of sequences of positive numbers whose mantissae are (or approach to be) distributed with respect to $\mu_{10}$, see e.g. \cite{CK, diaconis}. More recently, several authors have provided examples of sequences of random variables whose mantissa distribution converges to $\mu_{10}$ \cite{CMS, MS, sharpe} or whose the sequence of mantissae is almost surely distributed with respect to $\mu_{10}$. For a wide panorama on Benford sequences, see the reference books \cite{BH, miller}.

It is well known that a sequence $(x_n)$ of positive numbers is Benford in base $10$ if and only if the sequence of its fractional parts $(\{\log_{10} x_n\})$ is uniformly distributed in $[0, 1)$. According to the Weyl's criterion (see e.g. \cite{KN}, p7), the sequence $(x_n)$ is Benford if and only if, for any $h\in \ZZ^*$, we have
\[\lim_{N\rightarrow\infty} \frac{1}{N}\sum_{n=1}^Ne^{2i\pi h \log_{10} x_n}=0. \] 
To define a deviation between a sequence and the Benford's law, the notion of discrepancy is introduced as follows. Let $u=(u_n)$ be a sequence of real numbers. The discrepancy modulo 1 of order $N$ of $u$, associated with the natural density, is defined as
\[{D}_N(u) = \sup_{0\leq a<b<1}\left| \frac{1}{N} \sum_{n=1}^N \ind{[a,b)}(\{u_n\}) - (b-a)  \right|.\] For more details on the discrepancy, see e.g. \cite{KN}, p100--131.  For a sequence $x=(x_n)$, if we set $x_n=10^{u_n}$, we write $\overset{\sim}{D}_N(x) = {D}_N(u)$. The quantity $\overset{\sim}{D}_N(x)$ deals with the deviation between $\mu_{10}$ and the distribution of the first $N$ terms of $(\mathcal{M}_{10}(x_n))$ since $\{\log_{10} x_n\} = \log_{10}(\mathcal{M}_{10}(x_n))$. Hence 
\[\overset{\sim}{D}_N(x)=\sup_{1\leq s<t<10}\left|  \frac{1}{N}\sum_{n=1}^N \ind{[s,t)}(\mathcal{M}_{10}(x_n)) - \mu_{10}([s,t)) \right|.\]  In particular, $x=(x_n)$ is Benford if and only if $\overset{\sim}{D}_N(x)$ converges to 0 as $N$ goes to infinity. Through misuse of language, we also say that $\overset{\sim}{D}_N(x)$ is the discrepancy of $x=(x_n)$.

In this paper, we consider the following problem. Let $(X_n)$ be a sequence of positive independent random variables. We say that $(X_n)$ is a.s. Benford if $\omega-\PP a.s.$ the sequence $(X_n(\omega))$ is Benford. As observed in \cite{EMS}, several deterministic sequences at a power $d$ tend to be Benford when the power $d$ is large enough. The aim of our paper is to provide general conditions on the distribution of the random sequence $X=(X_n)$ to ensure that $X^{(d)}=(X_n^{d_n})$ is a.s. Benford for any sequence of positive numbers $(d_n)$ such that $d_n$ converges to infinity at a rate at most polynomial. 

First, we give some notation. In what follows, the function $\log$ denotes the natural logarithm.  For any functions $f$, $g$, we write $g(x)\underset{x\rightarrow\infty}{\sim}f(x)$ if and only if $\frac{g(x)}{f(x)}\conv[x]{\infty}1$. Moreover, we write $g(x)=O(f(x))$  if and only if there exists a positive number $M$ and a real number $x_0$ such that $|g(x)|\leq M|f(x)|$ for any $x\geq x_0$. 

We are now prepared to state our first theorem, which provides an upper bound for the discrepancy. 
 
\begin{Th}
\label{Prop:main2}
Let $(d_n)$ be a (deterministic) sequence of positive numbers such that $d_n=O\left(n^\theta\right)$ for some $\theta\geq 0$. Let $X=(X_n)$ be a sequence of positive independent random variables satisfying the   following two conditions:
\begin{enumerate}[(i)]
\item \label{cond1} there exists $\alpha>0$ such that $\sum_{n=1}^\infty \PPP{|\log X_n|>n^\alpha} < \infty;$
\item \label{cond2} there exists a sequence of nonnegative numbers $(r_n)$, with $r_n = O(n^{-\beta})$ for some $\beta>0$, and their exist four constants $c_1, c_2, \gamma, \delta>0$, such that for $n$ large enough and for each $h\in \NN^*$, we have 
\begin{equation}
\label{eq:cond2}
 \left|\EEE{e^{2i\pi h \log X_n}}\right|\leq c_1h^{-\gamma} + c_2h^\delta r_n.
 \end{equation}
\end{enumerate}
Then 
there exist an integrable random variable $C_0$ and a constant $c_0$ such that, for any $N\geq 1$, we have $\omega-\PP a.s. $ 
\begin{equation*}\overset{\sim}{D}_N(X^{(d)}(\omega))\leq C_0(\omega)\cdot (\log N)^2\cdot N^{-\frac{1}{2}} + c_0\left( \frac{1}{N}\sum_{n=1}^N(d_n)^{-\gamma}   + (\log N)^{\frac{1}{\delta +1}}\cdot N^{-\frac{\min\{\beta-\delta\theta, 1\}}{\delta + 1}}  \right),  
\end{equation*}
where $X^{(d)}(\omega) = (X_n^{d_n}(\omega))$. 
\end{Th}
The above theorem is obvious if the upper bound does not converge to 0. However, if $\delta\theta<\beta$, it provides a non-trivial estimate for the discrepancy when $d_n$ goes to infinity at a rate at most polynomial. As a consequence, we obtain the following result.

\begin{Cor}
\label{Cor1}
Let $(d_n)$ be such that $d_n=O\left(n^\theta\right)$ for some $\theta> 0$ and $d_n\conv[n]{\infty}\infty$. Assume that $X=(X_n)$  satisfies the assumptions \eqref{cond1} and \eqref{cond2} for some $\alpha, \beta, \gamma, \delta>0$, with $\delta\theta<\beta$. Then $\overset{\sim}{D}_N(X^{d}(\omega))$ converges $\omega-\PP$ a.s. to 0,  at a rate of convergence provided in Theorem \ref{Prop:main2}.  In particular, the sequence $(X_n^d(\omega))$ is a.s. Benford. 
\end{Cor}
In particular, if $X=(X_n)$  and $(d_n)$ satisfy the assumptions of Corollary \ref{Cor1}, with the more restrictive condition $d_n=O\left( n^\sigma \right)$ for each $\sigma>0$, then the discrepancy of $X^{(d)}(\omega)$ can be bounded as follows:
\[\sup_{1\leq s<t<10}\left|  \frac{1}{N}\sum_{n=1}^N \ind{[s,t)}(\mathcal{M}_{10}(X_n^{d_n}(\omega))) - \mu_{10}([s,t)) \right|\leq C(\omega)\cdot \frac{1}{N}\sum_{n=1}^Nd_n^{-\gamma}.\] 
 It is rather surprising that $X^{(d)}(\omega)$ is a.s. Benford for a sequence $d=(d_n)$ which converges arbitrarily slowly to infinity. On the opposite, it appears that for several classes of (deterministic) sequences $(x_n)$, the sequence $(x_n^{d_n})$ is Benford, when $(d_n)$ converges to infinity at a rate \emph{at less} polynomial (see e.g. Theorem 2 in \cite{MS2}). As a second consequence of Theorem \ref{Prop:main2}, the following corollary deals with the case where the sequence $(d_n)$ is constant.

\begin{Cor}
\label{Cor0}
Let $d_n=d$ for each $n\geq 1$ and let $X=(X_n)$ be such that the assumptions \eqref{cond1} and \eqref{cond2} hold for some $\alpha, \beta, \gamma, \delta>0$. Then there exist an integrable random variable $C_0(\omega)$ and a constant $c_0$ such that, for any $N\geq 1$, we have $\omega-\PP a.s. $ 
\[\overset{\sim}{D}_N(X^{d}(\omega)) \leq  C_0(\omega)\cdot (\log N)^2\cdot N^{-\frac{1}{2}} + c_0\left(d^{-\gamma} + (\log N)^{\frac{1}{\delta +1}}\cdot N^{-\frac{\min\{\beta, 1\}}{\delta + 1}}  \right),\]
where $X^{d}(\omega) = (X_n^d(\omega))$. 
\end{Cor}
In particular, as $d$ goes to infinity, the sequence $X^{d}=(X_n^d)$ tends to be a.s. Benford in the sense that its discrepancy converges to 0 as $d,N\rightarrow\infty$.  In a different context, such a convergence was already observed in Theorem 1 in \cite{EMS}, in which it is stated that two (deterministic)  sequences at a large power tend to be Benford. 

The assumption \eqref{cond1} of Theorem \ref{Prop:main2} is few restrictive. Indeed, thanks to the Markov's inequality, such a condition is satisfied when $\EEE{X_n}$ and $\EEE{X_n^{-1}}$ are negligible compared to $n^{-1-\epsilon}e^{n^\alpha}$ for some $\alpha, \epsilon>0$. The assumption \eqref{cond2}  of Theorem \ref{Prop:main2}  is in a way classical and is discussed in Remark \ref{Rk:assumption}.

Our paper is organized as follows. In Section \ref{sec:proof}, we prove Theorem \ref{Prop:main2}. This result is    illustrated through several examples of standard distributions in Section \ref{sec:examples}. These examples deal with discrete and continuous random variables respectively. In the rest of the paper, we denote by $c$ a generic constant which is independent of $\omega$, $N$ and $(d_n)$, but which may depend on other quantities.

\section{Proof of Theorem \ref{Prop:main2}}
\label{sec:proof}
To prove Theorem \ref{Prop:main2}, we apply two well-known inequalities. The first one deals with the discrepancy and is referred to as the Erd\"{o}s-Tur\'{a}n inequality (see e.g. \cite{RT}).

\begin{Th}(Erd\"{o}s-Tur\'{a}n inequality)
Let $x=(x_n)$ be a sequence of real numbers and let $N\geq 1$. Then, for every integer $H\geq 1$, we have
\[\overset{\sim}{D}_N(x) \leq \frac{1}{H+1} + \sum_{h=1}^H\frac{1}{h}\left| \frac{1}{N} \sum_{n=1}^Ne^{2i\pi h\log_{10} x_n}  \right|.\]
\end{Th}

The second inequality which we apply gives a deviation beween a sum of unit random complex numbers and the expectation of this sum. Such a result is due to Cohen and Cuny (Theorem 4.10 in \cite{CC}) and is re-written in our context. 

\begin{Th}(Cohen \& Cuny, 2006)
\label{Th:CohenCuny}
Let $(Y_n)$ be a sequence of independent random variables, with values in $\RR$. Assume that there exists $\eta>0$, such that $\sum_{n=1}^\infty \PPP{|Y_n|>n^\eta}<\infty$. Let $(a_n)$ be a sequence of complex numbers. Then there exist universal constants $\epsilon>0$ and $C>0$, such that 
\[\EEE{\sup_{N>K\geq 1} \sup_{T\geq 1}  \exp\left(\epsilon \cdot \frac{\max_{|t|\leq T} \left| \sum_{n=K+1}^N  a_n \left(e^{2i\pi t Y_n} - \EEE{e^{2i\pi t Y_n}}\right) \right|^2   }{\log(1+T)\log(1+N^\eta) \sum_{n=K+1}^N  |a_n|^2 } \right) } \leq C.\] 
\end{Th}

In the rest of the paper, with a slight abuse of notation, we omit the dependence in $\omega$, e.g. we write $\overset{\sim}{D}_N(X^{(d)})$ instead of $\overset{\sim}{D}_N(X^{(d)}(\omega))$. We are now prepared to prove our first theorem. 
\begin{prooft}{Theorem \ref{Prop:main2}}
  According to the Erd\"{o}s-Tur\'{a}n inequality, we have for any $H\geq 1$,
\[\overset{\sim}{D}_N(X^{(d)}) \leq \frac{1}{H+1} + \sum_{h=1}^H\frac{1}{h} \left| \frac{1}{N}\sum_{n=1}^N e^{2i\pi h\log_{10} X_n^{d_n}}   \right|.\]
Hence,
\begin{multline}
\label{eq:erdosturan}
\overset{\sim}{D}_N(X^{(d)}) \leq \frac{1}{H+1} + \sum_{h=1}^H\frac{1}{h} \left| \frac{1}{N}\sum_{n=1}^N \EEE{e^{2i\pi h\log_{10} X_n^{d_n}}}   \right|\\
 + \sum_{h=1}^H\frac{1}{h} \left| \frac{1}{N}\sum_{n=1}^N \left( e^{2i\pi h\log_{10} X_n^{d_n}} - \EEE{e^{2i\pi h\log_{10} X_n^{d_n}}}\right)  \right|.
\end{multline}

First, we provide an upper bound for the term on the bottom. To do it, we take $a_n=1$, $Y_n=\log_{10} X_n^{d_n}$ and $K=1$. Since $d_n=O(n^\theta)$, we obtain for $n$ large enough that $\PPP{|Y_n|>n^\eta}\leq \PPP{|\log X_n|>n^\alpha}$ with $\eta>\alpha+\theta$. Hence, according to the assumption \eqref{cond1}, we have $\sum_{n=1}^\infty \PPP{|Y_n|>n^\eta}<\infty$.  It follows from Theorem \ref{Th:CohenCuny} that 
\[\EEE{\sup_{N>1} \sup_{T\geq 1} \max_{|t|\leq T}  \frac{\left| \sum_{n=2}^N \left(e^{2i\pi t \log_{10} X_n^{d_n}} - \EEE{e^{2i\pi t \log_{10} X_n^{d_n}}}\right) \right|^2}{\log(1+T) \log(1+N^\eta)(N-1)}     } \leq C.\]
In particular, there exists an integrable random variable $c(\omega)$ such that, for any $N\geq 2$, $T\geq 1$, $|t|\leq T$ we have $\omega-\PP a.s. $ 
\[\left|\frac{1}{N} \sum_{n=1}^N \left(e^{2i\pi t \log_{10} X_n^{d_n}} - \EEE{e^{2i\pi t \log_{10} X_n^{d_n}}}\right)\right| \leq c(\omega) \cdot \sqrt{\log(1+T)}\cdot \sqrt{\frac{\log(1+N^\eta)}{N}}.  \] Notice that we have considered a sum over $n=1,\ldots, N$ and not over $n=2,\ldots, N$ in the above equation because $\left|  e^{2i\pi t \log_{10} X_1} - \EEE{e^{2i\pi t \log_{10} X_1}} \right|\leq 2$. By taking $T=H$ and $t=h$, we obtain for any $N\geq 1, H\geq 1$ that
\begin{multline}\label{eq:maj1} \sum_{h=1}^H\frac{1}{h} \left|\frac{1}{N} \sum_{n=1}^N \left(e^{2i\pi h \log_{10} X_n^{d_n}} - \EEE{e^{2i\pi h \log_{10} X_n^{d_n}}}\right)\right| \\
\begin{split}
 & \leq c(\omega)\sum_{h=1}^H\frac{1}{h}\sqrt{\log(1+H)}\cdot \sqrt{\frac{\log(1+N^\eta)}{N}}\\
&  \leq c'(\omega)\log H \sqrt{\log(1+H)}\cdot \sqrt{\frac{\log(1+N^\eta)}{N}}.  
\end{split}  \end{multline}

Secondly, we provide an upper bound for the second term in the right-hand side in \eqref{eq:erdosturan}. To do it, let $N_0$ be such that the inequality \eqref{eq:cond2} holds for each $N\geq N_0$. Then
\[\left|\frac{1}{N}\sum_{n=1}^N \EEE{e^{2i\pi h\log_{10} X_n^{d_n}}}   \right| \leq \frac{1}{N}\sum_{n=1}^{N_0} \left|\EEE{e^{2i\pi hd_n\log_{10} X_n}}   \right| + \frac{1}{N}\sum_{n=N_0+1}^N \left|\EEE{e^{2i\pi hd_n\log_{10} X_n}}   \right|.\]  Bounding $\left|\EEE{e^{2i\pi hd_n\log_{10} X_n}}   \right| $ by 1 in the first sum and applying the inequality \eqref{eq:cond2} in the second sum for the right-hand side, we get
\begin{equation*}
\left|\frac{1}{N}\sum_{n=1}^N \EEE{e^{2i\pi h\log_{10} X^{d_n}_n}}   \right| \leq \frac{N_0}{N} + c_1\cdot \frac{1}{N}\sum_{n=1}^N   \left(\frac{hd_n}{\log(10)}  \right)^{-\gamma} + c_2\cdot   \frac{1}{N}\sum_{n=1}^N  \left(\frac{hd_n}{\log(10)} \right)^\delta r_n.
\end{equation*}
Besides, $\sum_{h=1}^H\frac{1}{h}\leq c\log H$, $\sum_{h=1}^H\frac{1}{h^{1+\gamma}}\leq c$ and $\sum_{h=1}^H\frac{1}{h^{1-\delta}}\leq cH^{\delta}$. This implies that 
\begin{equation*}
\sum_{h=1}^H\frac{1}{h}\left|\frac{1}{N}\sum_{n=1}^N \EEE{e^{2i\pi h\log_{10} X_n^{d_n}}}   \right| \leq c\cdot \left(\frac{\log H}{N} + \frac{1}{N}\sum_{n=1}^N (d_n)^{-\gamma} +  \frac{1}{N}\sum_{n=1}^N(d_n)^\delta r_n\cdot H^\delta \right).
\end{equation*}
Since $d_n=O\left(n^\theta \right)$ and $r_n=O\left( n^{-\beta} \right)$, we have $\frac{1}{N}\sum_{n=1}^N(d_n)^\delta r_n \leq c\cdot \log N \cdot N^{-1}$ if $\beta-\delta\theta=1$ and  $\frac{1}{N}\sum_{n=1}^N(d_n)^\delta r_n \leq c\cdot  N^{-\min\{\beta - \delta\theta, 1\}}$ otherwise. This together with \eqref{eq:erdosturan} and  \eqref{eq:maj1} implies that 
\begin{multline*}\overset{\sim}{D}_N(X^{(d)})\leq \frac{1}{H+1} + c''(\omega)\cdot \log H\cdot \sqrt{\log(1+H)}\cdot \sqrt{\frac{\log(1+N^\eta)}{N}}\\ 
+ c\cdot \left( \frac{1}{N}\sum_{n=1}^N (d_n)^{-\gamma} +  \log N \cdot N^{-\min\{\beta - \delta\theta, 1\}}\cdot H^\delta \right).
\end{multline*}
Optimizing the right-hand side over $H\geq 1$, we conclude the proof of Theorem \ref{Prop:main2} by taking
\[H=\left\lfloor  (\log N)^{-\frac{1}{\delta +1}}\cdot  N^{\frac{\min\{\beta-\delta\theta, 1\}}{\delta+1}}  \right\rfloor +1.\]
\end{prooft}

\begin{Rk}
\label{Rk:assumption}
The assumption given in Equation \eqref{eq:cond2} has been chosen in such a way that it holds when $X_n$ follows the (discrete) uniform distribution on $\{1,\ldots, n\}$. Indeed, in this case, we have
\[ \left| \EEE{e^{2i\pi h\log X_n}}  \right|  = \left| \frac{1}{n} \sum_{k=1}^n e^{2i\pi h\log k}  \right| \leq \frac{1}{\sqrt{n}} + \frac{1}{n}\left| \sum_{k=\lfloor \sqrt{n} \rfloor + 1}^n e^{2i\pi h\log k}  \right| , \]
According to the Van der Corput's theorem (see e.g. \cite{KN}, p17), this shows that
\[\left| \EEE{e^{2i\pi h\log X_n}}  \right|  \leq \frac{8}{\sqrt{h}} + \frac{1+4\sqrt{h}}{\sqrt{n}}  + \frac{6}{n} + \frac{3h}{n\sqrt{n}}.\]
 In particular, this satisfies Equation \eqref{eq:cond2} with $\gamma=\frac{1}{2}$, $\delta=1$ and $r_n=\frac{1}{\sqrt{n}}$. However, our assumption \eqref{cond2} and our assumption on the independence of the random variables $X_n$ remain restrictive. We hope, in a future paper, to extent Theorem \ref{Prop:main2} with more general conditions. 
\end{Rk}

\begin{Rk}
\label{Rk:cuny}
The main tool to derive the rate of the discrepancy  is contained in Theorem \ref{Th:CohenCuny}. Besides, as a consequence of Corollary \ref{Cor0},   we deduce that  $\omega-\PP a.s. $  \begin{equation}
\label{eq:corollary} 
 \lim_{d\rightarrow\infty}\limsup_{N\rightarrow\infty}\overset{\sim}{D}_N(X^d)=0.
 \end{equation} In particular, when $d$ is large, the sequence $X^d=(X_n^d)$ tends to be a Benford sequence. However,  Theorem \ref{Th:CohenCuny} is not necessary to derive Equation \eqref{eq:corollary} because the latter can be proved directly by standard arguments. Indeed, it follows from the law of large numbers (for independent non-stationary random variables) and the Erd\"{o}s-Tur\'{a}n inequality that for all fixed $H\geq 1$, 
\[\limsup_{N\rightarrow \infty} \overset{\sim}{D}_N(X^d) \leq \frac{1}{H+1} + \sum_{h=1}^H\frac{1}{h}\limsup_{N\rightarrow\infty} \frac{1}{N} \left|\sum_{n=1}^N  \EEE{e^{2i\pi hd\log_{10} X_n}} \right|.\]
Besides, according to \eqref{eq:cond2}, we know that  \[\lim_{d\rightarrow\infty}\limsup_{N\rightarrow\infty} \frac{1}{N} \left|\sum_{n=1}^N  \EEE{e^{2i\pi hd\log_{10} X_n}} \right| = 0.\] Hence, by taking $H\rightarrow\infty$, this proves that $\lim_{d\rightarrow\infty}\limsup_{N\rightarrow\infty}\overset{\sim}{D}_N(X^d)=0$. However, the main contribution of our paper is to provide an explicit rate of convergence for the discrepancy of $X^d$ as $d$ goes to infinity. 
\end{Rk}


\section{Examples}
\label{sec:examples}
In this section, we give several examples of sequences of random variables satisfying the assumptions \eqref{cond1} and \eqref{cond2} of Theorem \ref{Prop:main2}. Our examples deal with discrete and continuous random variables respectively.

\subsection{Discrete random variables}
The following proposition provides sufficient conditions for discrete random variables to ensure that the  assumption \eqref{cond2} of Theorem \ref{Prop:main2} is satisfied for $\gamma=\delta=1$. 
\begin{Prop}
\label{Prop:discrete}
Let $(X_n)$ be a sequence of random variables with finite expectation and such that $X_n\geq 1$ a.s..  Assume that there exists a sequence of modes  $(m_n)$ such that the sequences  $(\PPP{X_n=k})_{k\leq m_n}$ and $(\PPP{X_n=k})_{k > m_n}$ are non-decreasing and non-increasing respectively. Moreover, assume that for some $\beta>0$ one of the two following cases is satisfied:
\begin{itemize}
\item \textbf{Case 1:} $m_n\cdot n^{-\beta}\conv[n]{\infty}\infty$  and $\sup_{n\geq 1}m_n\PPP{X_n=m_n}<\infty$;
\item \textbf{Case 2:} $\sup_{n\geq 1}m_n<\infty$, $\PPP{X_n=m_n} = O\left(n^{-\beta} \right)$ and $\EEE{\frac{1}{X_n}} = O\left( n^{-\beta} \right)$.
\end{itemize}
Then for $n$ large enough and for each $h\geq 1$, we have:
\[\left| \EEE{e^{2i\pi h\log X_n}} \right| \leq c_1h^{-1} + c_2hn^{-\beta}\] where $c_1,c_2$ are two constants. 
\end{Prop}

\begin{prooft}{Proposition \ref{Prop:discrete}}
First, we provide a generic upper bound for $\EEE{e^{2i\pi h\log X_n}}$ which is independent of the two above cases. Then we deduce a specific upper bound for this expectation which depends this time on the case which is considered. 

To do it,  we write $ \EEE{e^{2i\pi h\log X_n}} = \lim_{N\rightarrow\infty}\sum_{k=1}^Ne^{2i\pi h\log k}\PPP{X_n=k}$. Let $N\geq 1$ be fixed. It follows from the Abel transformation that
\begin{multline*}\sum_{k=1}^Ne^{2i\pi h\log k}\PPP{X_n=k} = \PPP{X_n=N+1}\sum_{j=1}^Ne^{2i\pi h\log j}\\
- \sum_{k=1}^{N}\sum_{j=1}^ke^{2i\pi h\log j}(\PPP{X_n=k+1}-\PPP{X_n=k}).
\end{multline*}
Since $\left|\PPP{X_n=N+1}\sum_{j=1}^Ne^{2i\pi h\log j}\right|\leq N\PPP{X_n=N+1}$ converges to 0 as $N$ goes to infinity (because $\EEE{X_n}<\infty$), it is enough prove that 
\[\left|  \sum_{k=1}^{N}\sum_{j=1}^ke^{2i\pi h\log j}(\PPP{X_n=k+1}-\PPP{X_n=k})  \right| \leq \frac{c_1}{h} + hc_2n^{-\beta},\] for some constants $c_1,c_2$. To do it, we apply the following lemma.

\begin{Le}
\label{Le:TAF}
For each $h\geq 1$, $k\geq 1$, we have
\[\left| \sum_{j=1}^k e^{2i\pi h\log j}  \right|  \leq \frac{k}{2\pi h} + 1 + \pi h\log k. \]
\end{Le}
\begin{prooft}{Lemma \ref{Le:TAF}}
First, we notice that
\[\sum_{j=1}^ke^{2i\pi h\log j} = k^{2i\pi h+1}R_k(f),\] where $R_k(f):=\sum_{j=0}^{k-1}\int_{\frac{j}{k}}^{\frac{j+1}{k}}f\left(\frac{j+1}{k} \right)\mathrm{d}t $ is the Riemann sum of the function $f: t\mapsto t^{2i\pi h}$ on $[0,1]$ with $n$ regular steps of length $n^{-1}$. Hence
\[\begin{split} \left|\sum_{j=1}^ke^{2i\pi h\log j}\right| & \leq  k\left| \int_0^1f(t)\mathrm{d}t  \right| +  k\left| R_k(f) - \int_0^1f(t)\mathrm{d}t  \right|\\
& \leq  \frac{k}{2\pi h} + k\left| R_k(f) - \int_0^1f(t)\mathrm{d}t  \right|,
\end{split}\]
where the second inequality comes from the fact that $\int_0^1f(t)\mathrm{d}t = \frac{1}{2i\pi h+1}$. Besides,
\[\begin{split}
\left| R_k(f) - \int_0^1f(t)\mathrm{d}t  \right| & = \left| \sum_{j=0}^{k-1} \int_{\frac{j}{k}}^{\frac{j+1}{k}} \left( f\left(\frac{j+1}{k}\right) - f(t) \right)\mathrm{d}t  \right|\\
& \leq \left| \int_0^{\frac{1}{k}} \left(f\left(\frac{1}{k}\right) - f(t)\right)\mathrm{d}t  \right| + \sum_{j=1}^{k-1}\int_{\frac{j}{k}}^{\frac{j+1}{k}}\left(\frac{j+1}{k}  - t \right)\cdot \frac{2\pi hk}{j}\mathrm{d}t,
\end{split}\]
where the last line is a consequence of the mean value inequality. Integrating the right-hand side over $t$, we get
\[\left| R_k(f) - \int_0^1f(t)\mathrm{d}t  \right| \leq \frac{1}{k} + 2\pi h\sum_{j=1}^{k-1}\frac{1}{2jk} \leq \frac{1}{k} + \pi h\cdot \frac{\log k}{k}.\] 
This concludes the proof of Lemma \ref{Le:TAF}. 
\end{prooft}
According to Lemma \ref{Le:TAF}, we have
\begin{multline*}
\left|  \sum_{k=1}^{N}\sum_{j=1}^ke^{2i\pi h\log j}(\PPP{X_n=k+1}-\PPP{X_n=k})  \right|\\
 \leq \sum_{k=1}^{N} \left(\frac{k}{2\pi h} + 1 + \pi h\log k \right) \left| \PPP{X_n=k+1} - \PPP{X_n=k}\right|. 
\end{multline*}
Since the sequences  $(\PPP{X_n=k})_{k\leq m_n}$ and $(\PPP{X_n=k})_{k\geq m_n}$ are non-decreasing and non-increasing respectively, we get
\begin{multline*}
\sum_{k=1}^{N} \left(\frac{k}{2\pi h} + 1 + \pi h\log k \right) \left| \PPP{X_n=k+1} - \PPP{X_n=k}  \right|\\
 = \sum_{k=1}^{m_n-1} \left(\frac{k}{2\pi h} +  \pi h\log k \right) \left( \PPP{X_n=k+1} - \PPP{X_n=k} \right)\\
 + \sum_{k=m_n}^{N} \left(\frac{k}{2\pi h} +  \pi h\log k \right) \left( \PPP{X_n=k} - \PPP{X_n=k+1}  \right)\\
 + 2\PPP{X_n=m_n}   - \PPP{X_n=N} - \PPP{X_n=1}. 
\end{multline*}
With standard computations, we get:
\begin{subequations}
\begin{equation*}
\sum_{k=1}^{m_n-1}  k\left( \PPP{X_n=k+1} - \PPP{X_n=k} \right)\leq m_n\PPP{X_n=m_n},
\end{equation*}
\begin{equation*}
\sum_{k=1}^{m_n-1}  \log k\left( \PPP{X_n=k+1} - \PPP{X_n=k} \right)\leq \log m_n \PPP{X_n=m_n},
\end{equation*}
\begin{equation*}
\sum_{k=m_n}^{N}  k\left( \PPP{X_n=k} - \PPP{X_n=k+1} \right)\leq m_n\PPP{X_n=m_n}+1,
\end{equation*}
\begin{multline*}
\sum_{k=m_n}^{N}  \log k\left( \PPP{X_n=k} - \PPP{X_n=k+1} \right)\leq \sum_{k=m_n}^{N-2}\log\left(1+\frac{1}{k} \right) \PPP{X_n=k+1}\\+\log m_n\PPP{X_n=m_n}.
\end{multline*}
\end{subequations}
Using the fact that $\log\left(1+\frac{1}{k}\right)\PPP{X_n=k+1}\leq \frac{1}{k}\PPP{X_n=k}$ for each $k\geq m_n$, we deduce that  
\begin{equation}
\label{eq:majdiscrete}
\sum_{k=1}^{N} \left(\frac{k}{2\pi h} + 1 +  \pi h\log k \right) \left| \PPP{X_n=k+1} - \PPP{X_n=k}  \right| \leq \frac{c_1}{h} + \pi hs_n,
\end{equation} where
\[c_1=\frac{1}{2\pi}\left(2\sup_{n\geq 1}m_n\PPP{X_n=m_n} + 1   \right)\] and
\[s_n =  2\log m_n \PPP{X_n=m_n}   + \sum_{k=m_n}^{N-1}\frac{1}{k}\PPP{X_n=k}  + 2\PPP{X_n=m_n}.\]
The inequality \eqref{eq:majdiscrete} is independent of the two cases considered in the assumptions of Proposition \ref{Prop:discrete}. Now, we deal with the terms $c_1$ and $s_n$ by discussing these two cases.
\begin{itemize}
\item Case 1: if $m_n\cdot n^{-\beta}\conv[n]{\infty}\infty$ for some $\beta>0$ and $\sup_{n\geq 1}m_n\PPP{X_n=m_n}<\infty$, 
we obtain that $c_1<\infty$. Moreover, $s_n = O\left( n^{-\beta} \right)$ since $\log m_n = O(m_n)$ and \[\sum_{k=m_n}^\infty \frac{1}{k}\PPP{X_n=k} \leq \sum_{k=m_n}^\infty\frac{1}{k}\PPP{X_n=m_n} \underset{n\rightarrow\infty}{\sim}\log m_n\cdot \PPP{X_n=m_n}.\]
\item Case 2: if $\sup_{n\geq 1}m_n<\infty$, $\PPP{X_n=m_n} = O\left(n^{-\beta}\right)$ and $\EEE{\frac{1}{X_n}} = O\left(n^{-\beta}\right)$ for some $\beta>0$, we also obtain that $c_1<\infty$ and $s_n = O\left(n^{-\beta}\right)$. 
\end{itemize}
This concludes the proof of Proposition \ref{Prop:discrete}. 
\end{prooft}
We give below three examples of sequences of random variables $X=(X_n)$ by checking the assumption \eqref{cond1} of Theorem \ref{Prop:main2} and one of the two cases of Proposition \ref{Prop:discrete}. According to Theorem \ref{Prop:main2} and Proposition \ref{Prop:discrete}, the discrepancy for each example can be bounded as follows:
\[\overset{\sim}{D}_N(X^{(d)}) \leq  C_0(\omega)\cdot  (\log N)^2\cdot N^{-\frac{1}{2}} + c_0\left( \frac{1}{N}\sum_{n=1}^N(d_n)^{-1} + (\log N)^{\frac{1}{2}}\cdot N^{-\frac{1}{2}\cdot\min\{\beta-\theta, 1\}}  \right).\]
In particular, if $(d_n)\rightarrow\infty$ with $d_n=O\left(n^\theta\right)$ and $\theta>\beta$, the sequence $X^{(d)} = (X_n^{d_n})$ is a.s. Benford.

\begin{Ex}
Assume that $X_n$ has a geometric distribution with parameter $p_n = O\left(n^{-\beta}\right)$. Here $m_n=1$, so that $\PPP{X_n=1}=p_n = O\left(n^{-\beta}\right)$. We also obtain the same order for 
$\EEE{\frac{1}{X_n}} = -\frac{p_n}{1-p_n}\cdot \log(1-p_n)$. In particular, the third conditions of Case 2 are satisfied. Besides, if $p_ne^{n^\alpha}n^{-\alpha'}\conv[n]{\infty}\infty$ for some $\alpha>0, \alpha'>1$, the assumption \eqref{cond1} holds since \[\sum_{n=1}^\infty \PPP{|\log X_n|>n^\alpha} \leq \sum_{n=1}^\infty \frac{1}{p_ne^{n^\alpha}}<\infty\] according to the Markov's inequality. 
\end{Ex}


\begin{Ex}
Let $X_n$ be a random variable with distribution $\PPP{X_n=k} = \frac{\alpha_n}{(n+k)^{1+\epsilon}}$, where $\alpha_n$ is the normalizing constant and $\epsilon>0$. In particular, we have 
\begin{equation}
\label{eq:exdiscrete2}
\epsilon n^\epsilon \leq \alpha_n\leq \epsilon (n+1)^\epsilon
\end{equation}
 since 
\begin{equation*}
\frac{1}{\int_n^{\infty} x^{-(1+\epsilon)}\mathrm{d}x  } \leq \alpha_n:=\frac{1}{\sum_{k=1}^\infty (n+k)^{-(1+\epsilon)}} \leq \frac{1}{\int_{n+1}^{\infty} x^{-(1+\epsilon)}\mathrm{d}x  }.
\end{equation*} Here $m_n=1$ and the third conditions of Case 2 are satisfied. Indeed, the first one is trivial and for the second one we have $\PPP{X_n=1} = O\left(n^{-(1+\epsilon)}\right)$. For the third condition, let $\beta <1$. According to \eqref{eq:exdiscrete2}, we have $\frac{1}{k}\cdot \frac{\alpha_n\cdot n^\beta}{(n+k)^{1+\epsilon}} \leq \frac{\epsilon}{k(k+1)^{1-\beta}}$. It follows from the dominated convergence theorem that 
\[\lim_{n\rightarrow\infty} n^\beta\cdot \EEE{\frac{1}{X_n}} = \sum_{k=1}^\infty\lim_{n\rightarrow\infty}\frac{1}{k}\cdot \frac{\alpha_n\cdot n^\beta}{(n+k)^{1+\epsilon}}=0.\] This checks the third condition of Case 2 for each $\beta<1$.  Besides, the assumption \eqref{cond1} holds since for each $n\geq 1$ and for each $\alpha>0$, we have
\[\PPP{|\log X_n |>n^\alpha} =  \sum_{k=\lfloor e^{n^\alpha} +1 \rfloor  }^\infty \frac{\alpha_n}{(n+k)^{1+\epsilon}}  \leq  \sum_{k=\lfloor e^{n^\alpha} +1 \rfloor  }^\infty \frac{\epsilon (n+1)^\epsilon}{(n+k)^{1+\epsilon}}   \underset{n\rightarrow\infty}{\sim} \frac{n^\epsilon}{e^{\epsilon n^{\alpha}}}.\]
\end{Ex}

\begin{Ex}
Assume that $X_n$ has a (discrete) uniform distribution in $\{a_n, \ldots, b_n\}$, with $a_n< b_n$,  $b_n\cdot n^{-\beta}\rightarrow\infty$ for some $\beta>0$, and $\limsup\frac{a_n}{b_n}< 1$. Here we take $m_n=b_n$. The two conditions of Case 1 are satisfied. Indeed, the first one holds because $b_n\cdot n^{-\beta}\rightarrow\infty$. The second one comes from the fact that $\limsup\frac{a_n}{b_n}< 1$ and $m_n\PPP{X_n=m_n} = \frac{b_n}{b_n-a_n+1}.$ Besides, a sufficient and few restrictive assumption on $b_n$ to ensure that the assumption \eqref{cond1} holds is:  $b_n=O(e^{n^\alpha})$ for some $\alpha >0$.   
Notice that if $\frac{a_n}{b_n}$ converges to 1, the random variables $X_n$ are asymptotically deterministic. It is not surprising that the property \eqref{Pty2} cannot hold in this context since there exist deterministic sequences such that, at any power $d$, the sequences are not Benford. 
\end{Ex}

\subsection{Continuous random variables}
Let $X=(X_n)$ be a sequence of random variables. We first state three properties which imply the  assumption \eqref{cond2} of Theorem \ref{Prop:main2} when they are simultaneously satisfied. 
\begin{enumerate}[(a)]
\item \label{Pty1} For any $n\geq 1$, the density $f_n$ of $X_n$ exists and is a piecewise absolutely continuous function. In what follows, we denote by $k_n$ the number of sub-domains of $f_n$ and by $I_{n,j}:=[a_{n,j}, b_{n,j}]$ the $j$-th sub-domain, with $a_{n,j}\leq b_{n,j}\leq a_{n,j+1}$ for each $1\leq j\leq k_n-1$. The $k_n$-th interval is of the form $I_{n,k_n}=[a_{n,k_n}, +\infty)$. In particular, $f_n$ is a.e. differentiable on  $\bigcup_{j=1}^{k_n} I_{n,j}$ and $f_n=0$ on the complement. 
\item \label{Pty2} $\limsup_{N\rightarrow\infty}\sum_{j=1}^{k_N}\sup_{x\in I_{N,j}} |xf_N(x)|<\infty$.
\item \label{Pty3} $\limsup_{N\rightarrow \infty} \sum_{j=1}^{k_N} \int_{I_{N,j}}|xf'_N(x)|\mathrm{d}x<\infty$. 
\end{enumerate}

Under the above assumptions, the following proposition ensures that the assumption \eqref{cond2} of Theorem \ref{Prop:main2} holds, with $\gamma=1$ and $a_n=0$ for each $n\geq 1$.

\begin{Prop}
\label{Prop:continuous}
If the properties hold \eqref{Pty1}, \eqref{Pty2} and \eqref{Pty3} hold simultaneously, then for $n$ large enough and for each $h\in \NN^*$, we have $\left| \EEE{e^{2i\pi h\log X_n}}  \right| \leq c_1h^{-1}$.
\end{Prop}

\begin{prooft}{Proposition \ref{Prop:continuous}}
It is enough to prove the following inequality:
\[\limsup_{N\rightarrow\infty}\sup_{h\in \NN^*}h\left|\EEE{e^{2i\pi h\log X_N}}\right|<\infty.\]
To do it, we assume without loss of generality that $k_n=1$ for each $n$, with $I_{n,j}=:I_n=[a_n,b_n]$. In particular, the density $f_n$ is  absolutely continuous on $[a_n,b_n]$ and equals 0 on the complement. This gives for any $N\geq 1, h\geq 1$
\[ \begin{split}  h\left|  \EEE{e^{2i\pi h\log X_N}} \right| & = h\left|  \int_{a_N}^{b_N} x^{2i\pi h}f_N(x)\mathrm{d}x    \right|\\
& = h\left|   \frac{1}{2i\pi h}\cdot \left(\left[x^{2i\pi h + 1}f_N(x) \right]_{a_N}^{b_N}  - \int_{a_N}^{b_N} x^{2i\pi h+1}f'_N(x)\mathrm{d}x  \right)\right|\\
& \leq \frac{1}{2\pi}\left(\sup_{x\in [a_N,b_N]}|xf_N(x)| + \int_{a_N}^{b_N}|xf'_N(x)|\mathrm{d}x \right).   
\end{split}\]
In particular, we have $\limsup_{N\rightarrow\infty}\sup_{h\in \NN^*}h\left|\EEE{e^{2i\pi h\log X_N}} \right|<\infty$ provided that the three above properties hold. 
\end{prooft}

Notice that if  $g_n$ denotes the density of $X_n^{-1}$, we can easily show that $g_n$ satisfies the above  assumptions if and only if the ones are satisfied by the density of $X_n$. This suggests that our assumptions are not very restrictive. We give below three examples of distributions of random variables which satisfy the assumption \eqref{cond1} of Theorem \ref{Prop:main2} and the three conditions \eqref{Pty1}, \eqref{Pty2} and \eqref{Pty3} of Proposition \ref{Prop:continuous}. According to Theorem \ref{Prop:main2} and Proposition \ref{Prop:continuous},  the discrepancy for each example can be bounded as follows:
\[\overset{\sim}{D}_N(X^{(d)}) \leq  C'_0(\omega)\cdot  (\log N)^2\cdot N^{-\frac{1}{2}} + c'_0\cdot \frac{1}{N}\sum_{n=1}^N(d_n)^{-1}.\]
To obtain the rate of the discrepancy, we have taken $\delta=1$ and $\beta\rightarrow\infty$. In particular, if $(d_n)\rightarrow\infty$ with $d_n=O\left(n^\theta\right)$ for some $\theta>0$, the sequence $X^{(d)} = (X_n^{d_n})$ is a.s. Benford.

\begin{Ex}
 If $X_n$ has an exponential distribution with parameter $\lambda_n>0$, the properties \eqref{Pty1}, \eqref{Pty2} and \eqref{Pty3} hold simultaneously, with $k_n=1$. Indeed, the first one is trivially satisfied and for the second and the third ones, we get: 
 \[\sup_{x\in \RR_+} |xf_n(x)| = e^{-1} \quad \text{and} \quad \int_{\RR_+}|xf'_n(x)|\mathrm{d}x = 1.\] Besides, for each $\alpha>0$, we have \[\PPP{|\log X_n|>n^\alpha} = e^{-\lambda_ne^{n^\alpha}} + (1-e^{-\lambda_ne^{-n^\alpha}}).\] Hence the assumption \eqref{cond1} is satisfied if there exists $\alpha'$ such that $\lambda_ne^{n^{\alpha'}}\conv[n]{\infty}\infty$ and $\lambda_ne^{-n^{\alpha'}}\conv[n]{\infty}0$. 
\end{Ex}

\begin{Ex}
 Assume that $X_n$ has a standard Fr\'echet distribution with parameter $\alpha_n>0$, i.e.  $\PPP{X_n\leq x} = e^{-x^{-\alpha_n}}$ if $x\geq 0$ and $\PPP{X_n\leq x} = 0$ otherwise. The property  \eqref{Pty1}  holds. Moreover, if $\inf_{n\geq 1}\alpha_n>0$ and $\sup_{n\geq 1}\alpha_n<\infty$, we can easily prove that the properties  \eqref{Pty2} and  \eqref{Pty3} are satisfied. Besides, the assumption \eqref{cond1} is also satisfied since for each $\alpha>0$, we have
 \[\PPP{|\log X_n|>n^\alpha} \underset{n\rightarrow\infty}{\sim}   e^{-\alpha_n\cdot n^\alpha} + e^{-e^{\alpha_n\cdot n^\alpha}},\] where the right-hand side is the term of a convergent series. 
\end{Ex}

\begin{Ex}
\label{Ex:uniform}
If $X_n$ has a (continuous) uniform distribution on $[a_n,b_n]$,  with $a_n< b_n$,  the properties  \eqref{Pty1} and  \eqref{Pty3} hold. Moreover, the property  \eqref{Pty2} is satisfied when $\limsup\frac{a_n}{b_n}< 1$. Besides, a sufficient and few restrictive assumption on $a_n,b_n$ to ensure that the assumption \eqref{cond1} holds is: $e^{-n^{\alpha}} = O(a_n)$ and $b_n=O(e^{n^\alpha})$ for some $\alpha >0$. Unsurprisingly, the assumptions on $b_n$ are very similar to those considered for a (discrete) uniform distribution.
\end{Ex}

\subsection{A numerical illustration}
In this section, we give a numerical illustration of a sequence of independent random variables $(X_n)$ such that  $(X_n^d)$ is almost a Benford sequence. For each $n$, the distribution of $X_n$ is assumed to be the (continuous) uniform distribution on $[1,n]$. This sequence satisfies the assumptions of Theorem \ref{Prop:main2} (see Example \ref{Ex:uniform}). In Table \ref{tab:uniform}, we provide the frequencies of the first significant digit of $X_1^d,\ldots, X_N^d$, with $N=1000$ and $d=2$. It appears that the distribution of frequencies of  $(X_n^d)$ is close to the Benford's law. 

\begin{table}[h!]
\begin{minipage}{7cm}\centering
\begin{tabular}{|c|c|c|}\hline
First digit & $(X_n^d)$ & Benford's law\\\hline    
1 & 0.308 & 0.306\\
2 & 0.204 & 0.184\\
3 & 0.096 & 0.116\\
4 & 0.116 & 0.106\\
5 &  0.084 & 0.082\\
6 & 0.068 & 0.055\\
7 & 0.060 & 0.050\\
8 & 0.028 & 0.053\\
9 &  0.036 & 0.048\\\hline
\end{tabular}\\
\end{minipage}
\begin{minipage}{7cm}\centering
 \includegraphics[height=5.0cm,width=5.5cm]{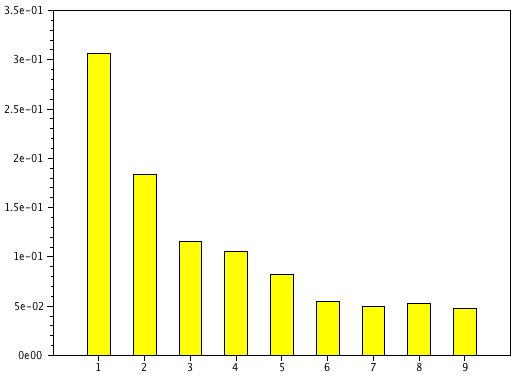} 
\end{minipage}
\caption{a simulation of the frequencies of the first significant digits of $X_1^d,\ldots, X_N^d$, where $X_n$ has a uniform distribution on $[1, n]$ for each $n\geq 1$, with $N=1000$ and $d=2$ (${\text{Scilab}}^\copyright$).}
\label{tab:uniform}
\end{table}


\end{document}